\documentclass[letterpaper, 10 pt, conference]{ieeeconf}%
\usepackage{graphicx}
\usepackage{epsfig}
\usepackage{amsmath}
\usepackage{amssymb}
\usepackage{bbm}
\usepackage{float}
\usepackage{tikz}
\usepackage{algorithm}
\usepackage{algpseudocode}
\usepackage{pgfplots}
\usepackage{gnuplottex}
\usepackage{mathrsfs}
\usepackage{subfigure}
\usepackage{booktabs}
\usepackage{multirow}
\usepackage{amsfonts}
\usepackage{wrapfig}
\usepackage{epstopdf}
\usepackage{bm}
\usetikzlibrary{plotmarks}
\usetikzlibrary{decorations.pathmorphing,shadows}

\IEEEoverridecommandlockouts
\allowdisplaybreaks
\providecommand{\U}[1]{\protect\rule{.1in}{.1in}}
\providecommand{\U}[1]{\protect\rule{.1in}{.1in}}
\makeatletter
\let\@ORGmakecaption\@makecaption
\long
\def\@makecaption#1#2{\@ORGmakecaption{#1}{#2}\vskip\belowcaptionskip\relax}
\makeatother
\def\@normalsize{\@setsize\normalsize{10pt}\xpt\@xpt
\abovedisplayskip 10pt plus2pt minus5pt\belowdisplayskip
\abovedisplayskip \abovedisplayshortskip \z@
plus3pt\belowdisplayshortskip 6pt plus3pt
minus3pt\let\@listi\@listI}
\newcounter{procount}
\newenvironment{proposition}[1][]{\refstepcounter{procount}\par   \textbf{Proposition~\theprocount. #1 } \rmfamily  }{\par}

\newcommand{\bsym}[1]{\boldsymbol{#1}}

\title{{\LARGE \textbf{Event excitation for event-driven control and optimization of multi-agent systems}}}
\author{ \parbox{3.5 in}{\centering Yasaman Khazaeni and Christos G. Cassandras\\
         Division of Systems Engineering\\ and Center for Information and Systems Engineering\\
         Boston University, MA 02446\\
         {\tt\small yas@bu.edu,cgc@bu.edu}}
         \thanks{The authors’ work is supported in part by NSF under grants CNS-
1239021, ECCS-1509084, and IIP-1430145, by AFOSR under grant FA9550-15-1-0471, by ONR under grant N00014-09-1-1051, and by the Cyprus Research Promotion
Foundation under Grant New Infrastructure Project/Strategic/0308/26.}
}

\begin{document}
\maketitle
\thispagestyle{empty}
\pagestyle{empty}
\begin{abstract}
We consider event-driven methods in a general framework for the control and optimization of multi-agent systems, viewing them as stochastic hybrid systems. Such systems often have feasible realizations in which the events needed to excite an on-line event-driven controller cannot occur, rendering the use of such controllers ineffective. We show that this commonly happens in environments which contain discrete points of interest which the agents must visit. To address this problem in event-driven gradient-based optimization problems, we propose a new metric for the objective function which creates a potential field guaranteeing that gradient values are non-zero when no events are present and which results in eventual event excitation. We apply this approach to the class of cooperative multi-agent data collection problems using the event-driven Infinitesimal Perturbation Analysis (IPA) methodology and include numerical examples illustrating its effectiveness.
\end{abstract}
\section{Introduction}

The modeling and analysis of dynamic systems has historically been founded on
the \emph{time-driven} paradigm provided by a theoretical framework based on
differential (or difference) equations: we postulate the existence of an
underlying ``clock'' and with every \textquotedblleft clock
tick\textquotedblright\ a state update is performed which synchronizes all
components of the system. As systems have become increasingly networked,
wireless, and distributed, the universal value of this paradigm has come to
question, since it may not be feasible to guarantee the synchronization of all
components of a distributed system, nor is it efficient to trigger actions
with every time step when such actions may be unnecessary. The
\emph{event-driven} paradigm offers an alternative to the modeling, control,
communication, and optimization of dynamic systems. The main idea in
event-driven methods is that actions affecting the system state need not be
taken at each clock tick. Instead, one can identify appropriate events that
trigger control actions. This approach includes the traditional time-driven
view if a clock-tick is considered a system \textquotedblleft
event\textquotedblright. Defining the right events is a crucial modeling step
and has to be carried out with a good understanding of the system dynamics.

The importance of event-driven behavior in dynamic systems was recognized in
the development of Discrete Event Systems (DES) and later Hybrid Systems (HS)
\cite{Cassandras2006}. More recently there have been significant advances in
applying event-driven methods (also referred to as ``event-based" and
``event-triggered") to classical feedback control systems; e.g., see
\cite{Heemels2008}, \cite{Tabuada10}, \cite{Trimpe2014}, as well as
\cite{Miskowicz2015event} and \cite{Cassandras2014event} and references
therein. Event-driven approaches are also attractive in receding horizon
control, where it is computationally inefficient to re-evaluate an optimal
control value over small time increments as opposed to event occurrences
defining appropriate planning horizons for the controller (e.g., see
\cite{KhazaeniCDC2014}). In distributed networked systems, event-driven
mechanisms have the advantage of significantly reducing communication among
networked components which cooperate to optimize a given objective.
Maintaining such cooperation normally requires frequent communication among
them; it was shown in \cite{Zhong2010asynchronous} that we can limit ourselves
to event-driven communication and still achieve optimization objectives while
drastically reducing communication costs (hence, prolonging the lifetime of a
wireless network), even when delays are present (as long as they are bounded).

Clearly, the premise of these methods is that the events involved are
observable so as to \textquotedblleft excite\textquotedblright\ the underlying
event-driven controller. However, it is not always obvious that these events
actually take place under every feasible control: it is possible that under
some control no such events are excited, in which case the controller may be
useless. In such cases, one can resort to artificial \textquotedblleft timeout
events\textquotedblright\ so as to eventually take actions, but this is
obviously inefficient. Moreover, in event-driven optimization mechanisms this
problem results in very slow convergence to an optimum or in an algorithm
failing to generate any improvement in the decision variables being updated.

In this work, we address this issue of event excitation in the context of
multi-agent systems. In this case, the events required are often defined by an
agent \textquotedblleft visiting\textquotedblright\ a region or a single point
in a mission space $S\subset\mathbb{R}^{2}$. Clearly, it is possible that such
events never occur for a large number of feasible agent trajectories. This is
a serious problem in trajectory planning and optimization tasks which are
common in multi-agent systems seeking to optimize different objectives
associated with these tasks, including coverage, persistent monitoring or
formation control \cite{Schwager2009decentralized}, \cite{Cassandras2013_2},
\cite{cao2011maintaining,Kwang2014,Yamaguchi1994,Desai1999,ji2007distributed,jiannan2013}%
. At the heart of this problem is the fact that objective functions for such
tasks rely on a non-zero reward (or cost) metric associated with a subset
$S^{+}\subset S$ of points, while all other points in $S$ have a reward (or
cost) which is zero since they are not \textquotedblleft points of
interest\textquotedblright\ in the mission space. We propose a novel metric
which allows \emph{all} points in $S$ to acquire generally non-zero reward (or
cost), thus ensuring that all events are ultimately excited. This leads to a
new method allowing us to apply event-based control and optimization to a
large class of multi-agent problems. We will illustrate the use of this method
by considering a general trajectory optimization problem in which
Infinitesimal Perturbation Analysis (IPA) \cite{Cassandras2006} is used as an
event-driven gradient estimation method to seek optimal trajectories for a
class of multi-agent problems where the agents must cooperatively visit a set
of target points to collect associated rewards (e.g., to collect data that are
buffered at these points.) This defines a family within the class of Traveling
Salesman Problems (TSPs) \cite{applegate2011} for which most solutions are
based on techniques typically seeking a shortest path in the underlying graph.
These methods have several drawbacks: $(i)$ they are generally combinatorially
complex, $(ii)$ they treat agents as particles (hence, not accounting for
limitations in motion dynamics which should not, for instance, allow an agent
to form a trajectory consisting of straight lines), and $(iii)$ they become
computationally infeasible as on-line methods in the presence of stochastic
effects such as random target rewards or failing agents. As an alternative we
seek solutions in terms of parameterized agent trajectories which can be
adjusted on line as a result of random effects and which are scalable, hence
computationally efficient, especially in problems with large numbers of
targets and/or agents. This approach was successfully used in \cite{Lin2015},
\cite{Khazaeni2015}.

In section \ref{Event-driven} we present the general framework for multi-agent
problems and address the event excitation issue. In section \ref{IPA} we
overview the event-driven IPA methodology and how it is applied to a general
hybrid system optimization problem. In section \ref{DH} we introduce a data
collection problem as an application of the general framework introduced in
section \ref{Event-driven} and will show simulation results of applying the
new methodology to this example in section \ref{numerical}.

\section{Event-Driven Optimization in Multi-Agent Systems}

\label{Event-driven} Multi-agent systems are commonly modeled as hybrid
systems with time-driven dynamics describing the motion of the agents or the
evolution of physical processes in a given environment, while event-driven
behavior characterizes events that may occur randomly (e.g., an agent failure)
or in accordance to control policies (e.g., an agent stopping to sense the
environment or to change directions). In some cases, the solution of a
multi-agent dynamic optimization problem is reduced to a policy that is
naturally parametric. As such, a multi-agent system can be studied with
parameterized controllers aiming to meet certain specifications or to optimize
a given performance metric. Moreover, in cases where such a dynamic
optimization problem cannot be shown to be reduced to a parametric policy,
using such a policy is still near-optimal or at least offers an alternative.

In order to build a general framework for multi-agent optimization problems,
assuming $S$ as the mission space, we introduce the function
$R(w):S\rightarrow\mathbb{R}$ as a \textquotedblleft property" of point $w\in
S$. For instance, $R(w)$ could be a weight that gives relative importance to
one point in $S$ compared to another. Setting $R(w)>0$ for only a finite
number of points implies that we limit ourselves to a finite set of points of
interest while the rest of $S$ has no significant value.

Assuming $F$ to be the set of all feasible agent states, We define
$P(w,s):S\times F\rightarrow\mathbb{R}$ to capture the cost/reward resulting
from how agents with state $s\in F$ interact with $w\in S$. For instance, in
coverage problems if an \textquotedblleft event\textquotedblright\ occurs at
$w$, then $P(w,s)$ is the probability of agents jointly detecting such events
based on the relative distance of each agent from $w$.

In general settings, the objective is to find the best state vector
$s_{1},\cdots,s_{N}$ so that $N$ agents achieve a maximal reward (minimal
cost) from interacting with the mission space $S$:
\begin{equation}
\min_{s \in F}J=\int_{S}P(w,s)R(w)dw \label{GeneralJ}%
\end{equation}

\begin{figure}
\centering
\includegraphics[width=2.4in,height=1.2in]{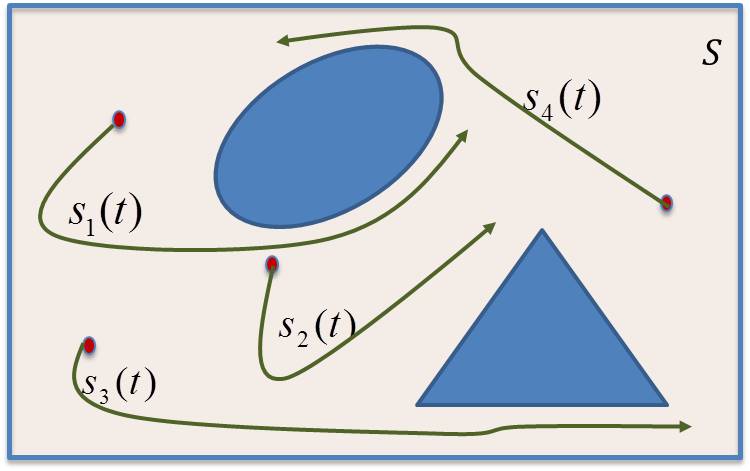}\newline%
\caption{Multi-agent system in a dynamic setting, blue areas are obstacles}%
\label{DynamicP}%
\end{figure}

This static problem can be extended to a dynamic version where the agents
determine optimal trajectories $s_{i}(t),$ $t\in\lbrack0,T]$, rather than
static states:
\begin{equation}
\min_{u(t)\in U}J=\int_{0}^{T}\int_{S}P(w,\mathbf{s}%
(u(t)))R(w,t)dwdt\label{GeneralJt}%
\end{equation}
subject to motion dynamics:
\begin{equation}
\dot{s}_{j}(t)=f_{j}(s_{j},u_{j},t),~j=1,\cdots,N\label{agentdynamic}%
\end{equation}
In Fig. \ref{DynamicP}, such a dynamic multi agent system is illustrated.
As an example, consensus problems are just a special case of \eqref{GeneralJ}.
Suppose that we consider a finite set of points $w\in S$ which coincide with
the agents states $s_{1},...,s_{N}$ (which are not necessarily their
locations). Then we can set $P(w,s)=\Vert s_{i}-s_{j}\Vert^{2}$ and,
therefore, replace the integral in \eqref{GeneralJ} by a sum. In this case,
$R(w)=R_{i}$ is just the weight that an agent carries in the consensus
algorithm. An optimum occurs when $\Vert s_{i}-s_{j}\Vert^{2}=0$ for all
$i,j$, i.e., all agents \textquotedblleft agree" and consensus is reached.
This is a special case because of the simplicity in $P(w,s)$ making the
problem convex so that a global optimum can be achieved, in contrast to most
problems we are interested in.

As for the formulation in \eqref{GeneralJt}, consider a trajectory planning
problem where $N$ mobile agents are tasked to visit $M$ stationary targets in
the mission space $S$. Target behavior is described through state variables
$x_{i}(t)$ which may model reward functions, the amount of data present at
$i$, or other problem-dependent target properties.
%
More formally, let $(\Omega,\mathcal{F},\mathcal{P})$ be an appropriately defined probability space
and $\omega\in\Omega$ a realization of the system where target dynamics are
subject to random effects:
\begin{equation}
\dot{x}_{i}(t)=g_{i}(x_{i}(t),\omega)\label{xdynamics}%
\end{equation}
$g_{i}(\cdot)$ is as such that $x_{i}(t)$ is monotonically increasing by $t$
and it resets to zero each time a target is completely emptied by an agent. In
the context of \eqref{GeneralJt}, we assume the $M$ targets are located at
points $w_{i},~i=1,\cdots,M$ and define
\begin{equation}
R(w,t)=\left\{
\begin{array}
[c]{ll}%
R(x_{i}(t),w) & \mbox{if }w\in C(w_{i})\\
0 & \mbox{otherwise}
\end{array}
\right.  \label{Rfunction_finite}%
\end{equation}
to be the value of point $w$, where $C(w_{i})$ is a compact 2-manifold in $\mathbb{R}^{2}$ containing $w_{i}$ which can be considered to be a region
defined by the sensing range of that target relative to agents (e.g., a disk
centered at $w_{i}$). Note that $R(w,t)$ is also a random variable defined on
the same probability space above. Given that only points $w\in C(w_{i})$ have
value for the agents, there is an infinite number of points $w\notin C(w_{i})$
such that $R(w,t)=0$ provided the following condition holds:

\textbf{Condition 1:} If $\exists i~\mbox{such that}~w\in C(w_{i})$ then
$w\notin C(w_{j})$ holds $\forall j\neq i$.

This condition is to ensure that two targets do not share any point $w$ in
their respective sensing ranges. Also it ensures that the set $\{C(w_{i}) ~|~
i=1:\cdots,M\}$ does not create a compact partitioning of the mission space
and there exist points $w$ which do not belong to any of the $C(w_{i})$.

Viewed as a stochastic hybrid system, we may define different modes depending
on the states of agents or targets and events that cause transitions between
these modes. Relative to a target $i$, any agent has at least two modes: being
at a point $w\in C(w_{i})$, i.e., visiting this target or not visiting it.
Within each mode, agent $j$'s dynamics, dictated by \eqref{agentdynamic}, and
target $i$'s dynamics in (\ref{xdynamics}) may vary. Accordingly, there are at
least two types of events in such a system: $(i)$ $\delta_{ij}^{0}$ events
occur when agent $j$ initiates a visit at target $i$, and $(ii)$ $\delta
_{ij}^{+}$ events occur when agent $j$ ends a visit at target $i$. Additional
event types may be included depending on the specifics of a problem, e.g.,
mode switches in the target dynamics or agents encountering obstacles.

\begin{figure}[ptb]
\centering
\includegraphics[width=2.7in,height=1.1in]{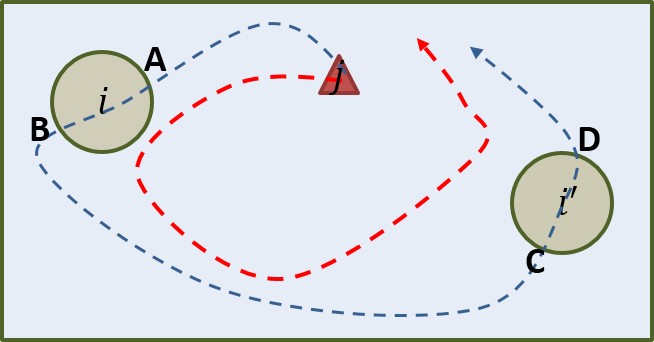}\newline%
\caption{Sample trajectories}%
\label{SampleEvents}%
\end{figure}

An example is shown in Fig. \ref{SampleEvents}, where target sensing ranges
are shown with green circles and agent trajectories are shown in dashed lines
starting at a base shown by a red triangle. In the blue trajectory, agent $1$
moves along the trajectory that passes through points $A\rightarrow
B\rightarrow C\rightarrow D$. It is easy to see that when passing through
points $A$ and $C$ we have $\delta_{i1}^{0}$ and $\delta_{i^{\prime}1}^{0}$
events, while passing through $B$ and $D$ we have $\delta_{i1}^{+}$ and
$\delta_{i^{\prime}1}^{+}$ events. The red trajectory is an example where none
of the events is excited. Suppose we consider an on-line trajectory adjustment
process in which the agent improves its trajectory based on its performance
measured through \eqref{Rfunction_finite}. In this case, $R(w,t)=0$ over all
$t$, as long as the agent keeps using the red trajectory, i.e., no event ever
occurs. Therefore, if an event-driven approach is used to control the
trajectory adjustment process, no action is ever triggered and the approach is
ineffective. In contrast, in the blue trajectory the controller can extract
useful information from every observed event; such information (e.g., a
gradient of $J$ with respect to controllable parameters as described in the
next section) can be used to adjust the current trajectory so as to improve
the objective function $J$ in (\ref{GeneralJ}) or (\ref{GeneralJt}).

Therefore, if we are to build an optimization framework for this class of
stochastic hybrid systems to allow the application of event-driven methods by
calculating a performance measure gradient, then a fundamental property
required is the occurrence of at least some events in a sample realization. In
particular, the IPA method \cite{Cassandras2010} is based on a single sample
realization of the system over which events are observed along with their
occurrence times and associated system states. Suppose that the trajectories
can be controlled through a set of parameters forming a vector $\theta$. Then,
IPA provides an unbiased estimate of the gradient of a performance metric
$J(\theta)$ with respect to $\theta$. This gradient is then used to improve
the trajectory and ultimately seek an optimal one when appropriate conditions hold.

As in the example of Fig. \ref{SampleEvents}, it is possible to encounter
trajectory realizations where no events occur in the system. In the above
example, this can easily happen if the trajectory does not pass through any
target. The existence of such undesirable trajectories is the direct
consequence of Condition 1. This lack of event excitation results in
event-based controllers being unsuitable.

\textbf{New Metric: } In order to overcome this issue we propose a new
definition for $R(w,t)$ in (\ref{Rfunction_finite}) as follows:
\begin{equation}
R(w,t)=\sum_{i=1}^{M}h_{i}(x_{i}(t),d_{i}(w)) \label{Rfunction}%
\end{equation}
where $w\in S$, $h_{i}(\cdot)$ is a function of the target's state $x_{i}(t)$
and $d_{i}(w)=\Vert w_{i}-w\Vert$. Note that, if $h_{i}(\cdot)$ is properly
defined, \eqref{Rfunction} yields $R(w,t)>0$ at all points.

While the exact form of $h_{i}(\cdot)$ depends on the problem, we impose the
condition that $h_{i}(\cdot)$ is monotonically decreasing in $d_{i}(w)$. We
can think of $h_{i}(\cdot)$ as a value function associated with point $w_{i}$.
Using the definition of $R(w,t)$, this value is spread out over all points
$w\in S$ rather than being concentrated at the single point $w_{i}$. This
creates a continuous potential field for the agents leading to a non-zero
gradient of the performance measure even when the trajectories do not excite
any events. This non-zero gradient will then induce trajectory adjustments
that naturally bring them toward ones with observable events.

Finally, recalling the definition in \eqref{GeneralJt}, we also define:
\begin{equation}
P(w,s)=\sum_{j=1}^{N}\Vert s_{j}(t)-w\Vert^{2} \label{Pfunction}%
\end{equation}
the total quadratic travel cost for agents to visit point $w$.

In Section IV, we will show how to apply $R(w,t)$ and $P(w,s)$ defined as
above in order to determine optimal agent trajectories for a class of
multi-agent problems of the form (\ref{GeneralJt}). First, however, we review
in the next section the event-driven IPA calculus which allows us to estimate
performance gradients with respect to controllable parameters.

\section{Event-Driven IPA Calculus}

\label{IPA} Let us fix a particular value of the parameter $\boldsymbol{\theta
} \in\Theta$ and study a resulting sample path of a general SHS. Over such a
sample path, let $\tau_{k}(\boldsymbol{\theta})$, $k = 1, 2, \cdots$ denote
the occurrence times of the discrete events in increasing order, and define
$\tau_{0}(\boldsymbol{\theta}) = 0$ for convenience. We will use the notation
$\tau_{k}$ instead of $\tau_{k}(\boldsymbol{\theta})$ when no confusion
arises. The continuous state is also generally a function of
$\boldsymbol{\theta}$, as well as of $t$, and is thus denoted by
$x(\boldsymbol{\theta}, t)$. Over an interval $[\tau_{k}(\boldsymbol{\theta}),
\tau_{k+1}(\boldsymbol{\theta}))$, the system is at some mode during which the
time-driven state satisfies $\dot x = f_{k}(x, \boldsymbol{\theta}, t)
\label{dynamic} $, in which $x$ is any of the continuous state variables of
the system and $\dot x$ denotes $\frac{\partial x}{\partial t}$ . Note that we
suppress the dependence of $f_{k}$ on the inputs $u \in U$ and $d \in D$ and
stress instead its dependence on the parameter $\boldsymbol{\theta}$ which may
generally affect either $u$ or $d$ or both. The purpose of perturbation
analysis is to study how changes in $\boldsymbol{\theta}$ influence the state
$x(\boldsymbol{\theta}, t)$ and the event times $\tau_{k}(\boldsymbol{\theta
})$ and, ultimately, how they influence interesting performance metrics that
are generally expressed in terms of these variables.

An event occurring at time $\tau_{k+1}(\boldsymbol{\theta})$ triggers a change
in the mode of the system, which may also result in new dynamics represented
by $f_{k+1}$. The event times ${\tau_{k}(\boldsymbol{\theta})}$ play an
important role in defining the interactions between the time-driven and
event-driven dynamics of the system.

Following the framework in \cite{Cassandras2010}, consider a general
performance function $J$ of the control parameter $\boldsymbol{\theta}$:
\begin{equation}
J(\boldsymbol{\theta};x(\boldsymbol{\theta},0),T)=E[L(\boldsymbol{\theta
};x(\boldsymbol{\theta},0),T)]
\end{equation}
where $L(\boldsymbol{\theta};x(\boldsymbol{\theta},0),T)$ is a sample function
of interest evaluated in the interval $[0,T]$ with initial conditions
$x(\theta,0)$. For simplicity, we write $J(\boldsymbol{\theta})$ and
$L(\boldsymbol{\theta})$. Suppose that there are $K$ events, with occurrence
times generally dependent on $\boldsymbol{\theta}$, during the time interval
$[0,T]$ and define $\tau_{0}=0$ and $\tau_{N+1}=T$. Let $L_{k}:\mathbb{R}^{n}%
\times\Theta\times\mathbb{R}^{+}\rightarrow\mathbb{R}$ be a function and
define $L(\boldsymbol{\theta})$ by
\begin{equation}
L(\boldsymbol{\theta})=\sum_{k=0}^{K}\int_{\tau_{k}}^{\tau_{k+1}}%
L_{k}(x,\theta,t)dt \label{samplefunction}%
\end{equation}
where we reiterate that $x=x(\boldsymbol{\theta},t)$ is a function of
$\boldsymbol{\theta}$ and $t$. We also point out that the restriction of the
definition of $J(\boldsymbol{\theta})$ to a finite horizon $T$ which is
independent of $\boldsymbol{\theta}$ is made merely for the sake of
simplicity. Returning to the stochastic setting, the ultimate goal of the
iterative process shown is to maximize $E_{\omega}[L(\boldsymbol{\theta
},\omega)]$, where we use $\omega$ to emphasize dependence on a sample path
$\omega$ of a SHS (clearly, this is reduced to $L(\boldsymbol{\theta})$ in the
deterministic case). Achieving such optimality is possible under standard
ergodicity conditions imposed on the underlying stochastic processes, as well
as the assumption that a single global optimum exists; otherwise, the
gradient-based approach is simply continuously attempting to improve the
observed performance $L(\boldsymbol{\theta},\omega)$. Thus, we are interested
in estimating the gradient
\begin{equation}
\frac{dJ(\boldsymbol{\theta})}{d\boldsymbol{\theta}}=\frac{dE_{\omega
}[L(\boldsymbol{\theta},\omega)]}{d\boldsymbol{\theta}} \label{Jgradient}%
\end{equation}
by evaluating $\frac{dL(\boldsymbol{\theta},\omega)}{d\theta}$ based on
directly observed data. We obtain $\boldsymbol{\theta}^{\ast}$ by optimizing
$J(\boldsymbol{\theta})$ through an iterative scheme of the form
\begin{equation}
\boldsymbol{\theta}_{n+1}=\boldsymbol{\theta}_{n}-\eta_{n}H_{n}%
(\boldsymbol{\theta}_{n};x(\boldsymbol{\theta},0),T,\omega_{n}),~n=0,1,\cdots
\label{gradientmethod}%
\end{equation}
where ${\eta_{n}}$ is a step size sequence and $H_{n}(\boldsymbol{\theta}%
_{n};x(\boldsymbol{\theta},0),T,\omega_{n})$ is the estimate of $\frac
{dJ(\boldsymbol{\theta})}{d\boldsymbol{\theta}}$ at $\boldsymbol{\theta
}=\boldsymbol{\theta}_{n}$. In using IPA, $H_{n}(\boldsymbol{\theta}%
_{n};x(\boldsymbol{\theta},0),T,\omega_{n})$ is the sample derivative
$\frac{dL(\boldsymbol{\theta},\omega)}{d\boldsymbol{\theta}}$ , which is an
unbiased estimate of $\frac{dJ(\boldsymbol{\theta})}{d\boldsymbol{\theta}}$ if
the condition (dropping the symbol $\omega$ for simplicity)%
\begin{equation}
E\big[\frac{dL(\boldsymbol{\theta})}{d\boldsymbol{\theta}}\big]=\frac{d
E[L(\boldsymbol{\theta})]}{d\boldsymbol{\theta}} =\frac{dJ(\boldsymbol{\theta
})}{d\boldsymbol{\theta}} \label{condition}%
\end{equation}
is satisfied, which turns out to be the case under mild technical conditions.
The conditions under which algorithms of the form \eqref{gradientmethod}
converge are well-known (e.g., see \cite{Kushner2003}). Moreover, in addition
to being unbiased, it can be shown that such gradient estimates are
independent of the probability laws of the stochastic processes involved and
require minimal information from the observed sample path. The process through
which IPA evaluates $\frac{dL(\boldsymbol{\theta})}{d\boldsymbol{\theta}}$ is
based on analyzing how changes in $\boldsymbol{\theta}$ influence the state
$x(\boldsymbol{\theta}, t)$ and the event times $\tau_{k}(\boldsymbol{\theta
})$. In turn, this provides information on how $L(\boldsymbol{\theta})$ is
affected, because it is generally expressed in terms of these variables. Given
$\boldsymbol{\theta} = [\theta_{1}, . . . , \theta_{l} ]^{T}$, we use the
Jacobian matrix notation:
\begin{equation}
x^{\prime}(\boldsymbol{\theta}, t) = \frac{\partial x(\boldsymbol{\theta},
t)}{\partial\boldsymbol{\theta}},~ {\tau_{k}}^{\prime}= \frac{\partial\tau
_{k}(\boldsymbol{\theta})}{\partial\boldsymbol{\theta}}, k = 1, \cdots, K
\end{equation}
for all state and event time derivatives. For simplicity of notation, we omit
$\theta$ from the arguments of the functions above unless it is essential to
stress this dependence. It is shown in \cite{Cassandras2010} that $x^{\prime
}(t)$ satisfies:
\begin{equation}
\frac{dx^{\prime}(t)}{dt} =\frac{\partial f_{k}(t)}{\partial x}x^{\prime
}(t)+\frac{\partial f_{k}(t)}{\partial\boldsymbol{\theta}} \label{IPAderx}%
\end{equation}
for $t \in[\tau_{k}(\theta), \tau_{k+1|}(\theta))$ with boundary condition
\begin{equation}
x^{\prime}(\tau_{k}^{+} ) = x^{\prime}(\tau_{k}^{-}) + [ f_{k - 1}(\tau
_{k}^{-} ) - f_{k}(\tau_{k}^{+})]\tau_{k}^{\prime} \label{IPAderjump}%
\end{equation}
for $k = 0, \cdots, K$. We note that whereas $x(t)$ is often continuous in
$t$, $x^{\prime}(t)$ may be discontinuous in $t$ at the event times $\tau_{k}%
$; hence, the left and right limits above are generally different. If $x(t)$
is not continuous in $t$ at $t = \tau_{k}(\boldsymbol{\theta})$, the value of
$x(\tau_{k}^{+})$ is determined by the reset function $r(q, q^{\prime}, x,
\nu, \delta)$ and
\begin{equation}
x^{\prime}(\tau_{k}^{+} ) =\frac{d r(q, q^{\prime}, x, \nu, \delta
)}{d\boldsymbol{\theta}}%
\end{equation}

Furthermore, once the initial condition $x^{\prime}(\tau_{k}^{+})$ is given,
the linearized state trajectory ${x^{\prime}(t)}$ can be computed in the
interval $t\in\lbrack\tau_{k}(\boldsymbol{\theta}),\tau_{k+1}%
(\boldsymbol{\theta}))$ by solving \eqref{IPAderx} to obtain
\begin{equation}
x^{\prime}(t)=e^{\int_{\tau_{k}}^{t}\frac{\partial f_{k}(u)}{\partial x}%
du}\Big[\int_{\tau_{k}}^{t}\frac{\partial f_{k}(v)}{\partial\boldsymbol{\theta
}}e^{-\int_{\tau_{k}}^{t}\frac{\partial f_{k}(u)}{\partial x}du}dv+\xi
_{k}\Big] \label{xprime}%
\end{equation}
with the constant $\xi_{k}$ determined from $x^{\prime}(\tau_{k}^{+})$. In
order to complete the evaluation of $x^{\prime}(\tau_{k}^{+})$ we need to also
determine $\tau_{k}^{\prime}$. If the event at $\tau_{k}(\theta)$ is exogenous
$\tau_{k}^{\prime}=0$ and if the event at $\tau_{k}(\theta)$ is endogenous:
\begin{equation}
\tau_{k}^{\prime}=-\Big[\frac{\partial g_{k}}{\partial x}f_{k}(\tau_{k}%
^{-})\Big]\big(\frac{\partial g_{k}}{\partial\boldsymbol{\theta}}%
+\frac{\partial g_{k}}{\partial x}x^{\prime}(\tau_{k}^{-}%
)\big) \label{taukprime}%
\end{equation}
where $g_{k}(x,\boldsymbol{\theta})=0$ and it is defined as long as
$\frac{\partial g_{k}}{\partial x}f_{k}(\tau_{k}^{+})\neq0$ (details may be
found in \cite{Cassandras2010}.)

The derivative evaluation process involves using the IPA calculus in order to
evaluate the IPA derivative $\frac{dL}{d\boldsymbol{\theta}}$. This is
accomplished by taking derivatives in \eqref{samplefunction} with respect to
$\boldsymbol{\theta}$:
\begin{equation}
\frac{dL(\boldsymbol{\theta})}{d\boldsymbol{\theta}}=\sum_{k=0}^{K}\frac
{d}{d\boldsymbol{\theta}}\int_{\tau_{k}}^{\tau_{k+1}}L_{k}%
(x,\boldsymbol{\theta},t)dt
\end{equation}
Applying the Leibnitz rule, we obtain, for every $k=0,\cdots,K,$
\begin{equation}%
\begin{split}
\frac{d}{d\boldsymbol{\theta}}  &  \int_{\tau_{k}}^{\tau_{k+1}}L_{k}%
(x,\boldsymbol{\theta},t)dt\\
&  =\int_{\tau_{k}}^{\tau_{k+1}}\Big[\frac{\partial L_{k}(x,\boldsymbol{\theta
},t)}{\partial x}x^{\prime}(t)+\frac{\partial L_{k}(x,\boldsymbol{\theta}%
,t)}{\partial\boldsymbol{\theta}}\Big]dt\\
&  +L_{k}(x(\tau_{k+1}),\boldsymbol{\theta},\tau_{k+1})\tau_{k+1}^{\prime
}-L_{k}(x(\tau_{k}),\boldsymbol{\theta},\tau_{k})\tau_{k}^{\prime
}\label{IPAsampleder}%
\end{split}
\end{equation}

In summary the three equations \eqref{IPAderjump}, \eqref{xprime} and
\eqref{taukprime} form the basis of the IPA calculus and allow us to calculate
the final derivative in \eqref{IPAsampleder}. In the next section IPA is
applied to a data collection problem in a multi-agent system.

\section{The Data Collection Problem}

\label{DH} We consider a class of multi-agent problems where the agents must
cooperatively visit a set of target points to collect associated rewards
(e.g., to collect data that are buffered at these points.). The mission space
is $S\subset\mathbb{R}^{2}$. This class of problems falls within the general
formulation introduced in \eqref{GeneralJt}. The state of the system is the
position of agent $j$ time $t$, $s_{j}(t)=[s_{j}^{x}(t),s_{j}^{y}(t)]$ and the
state of the target $i$, $x_{i}(t)$. The agent's dynamics \eqref{agentdynamic}
follow a single integrator:
\begin{equation}
\dot{s}_{j}^{x}(t)=u_{j}(t)\cos\theta_{j}(t),\quad\ \ \dot{s}_{j}^{y}%
(t)=u_{j}(t)\sin\theta_{j}(t)\label{agentdynamics}%
\end{equation}
where $u_{j}(t)$ is the scalar speed of the agent (normalized so that $0\leq
u_{j}(t)\leq1$) and $\theta_{j}(t)$ is the angle relative to the positive
direction, $0\leq\theta_{j}(t)<2\pi$. Thus, we assume that each agent controls
its speed and heading.

We assume the state of the target $x_{i}(t)$ represents the amount of data
that is currently available at target $i$ (this can be modified to different
state interpretations). The dynamics of $x_{i}(t)$ in (\ref{xdynamics}) for
this problem are:
\begin{equation}
\resizebox{0.99 \columnwidth}{!}{$ \dot{x}_{i}(t)=\left\{ \begin{array}{ll} 0\qquad \mbox{if }x_{i}(t)=0\mbox{ and }\sigma _{i}(t)\leq \mu _{ij}p(s_j(t),w_i) & \\ \sigma _{i}(t)-\mu _{ij}p(s_j(t),w_i) \qquad\qquad \qquad \mbox{otherwise} & \end{array}\right. $}\label{Xdot}%
\end{equation}
i.e., we model the data at the target as satisfying simple flow dynamics with
an exogenous (generally stochastic) inflow $\sigma_{i}(t)$ and a controllable
rate with which an agent empties the data queue given by $\mu_{ij}%
p(s_{j}(t),w_{i})$. For brevity we set $p(s_{j}(t),w_{i})=p_{ij}(t)$ which is
the normalized data collection rate from target $i$ by agent $j$ and $\mu
_{ij}$ is a nominal rate corresponding to target $i$ and agent $j$.

Assuming $M$ targets are located at $w_{i}\in S,$ $i=1,\dots,M,$ and have a
finite range of $r_{i}$, then agent $j$ can collect data from $w_{i}$ only if
$d_{ij}(t)=\Vert w_{i}-s_{j}(t)\Vert\leq r_{i}$. We then assume that:
$(\mathbf{A1})$ $p_{ij}(t)\in\lbrack0,1]$ is monotonically non-increasing in
the value of $d_{ij}(t)=\Vert w_{i}-s_{j}(t)\Vert$, and $(\mathbf{A2})$ it
satisfies $p_{ij}(t)=0$ if $d_{ij}(t)>r_{i}$. Thus, $p_{ij}(t)$ can model
communication power constraints which depend on the distance between a data
source and an agent equipped with a receiver (similar to the model used in
\cite{Ny2008}) or sensing range constraints if an agent collects data using
on-board sensors. For simplicity, we will also assume that: $(\mathbf{A3})$
$p_{ij}(t)$ is continuous in $d_{ij}(t)$ and $(\mathbf{A4})$ only one agent at
a time is connected to a target $i$ even if there are other agents $l$ with
$p_{il}(t)>0$; this is not the only possible model, but we adopt it based on
the premise that simultaneous downloading of packets from a common source
creates problems of proper data reconstruction. This means that $j$ in
\eqref{Xdot} is the index of the agent that is connected to target $i$ at time
$t$.

The dynamics of $x_{i}(t)$ in \eqref{Xdot} results in two new event types
added to what was defined earlier, $(i)$ $\xi_{i}^{0}$ events occur when
$x_{i}(t)$ reaches zero, and $(ii)$ $\xi_{i}^{+}$ events occur when $x_{i}(t)$ leaves zero.

The performance measure is the total content of data left at targets at the
end of a finite mission time $T$. Thus, we define $J_{1}(t)$ to be the
following (recalling that $\{\sigma_{i}(t)\}$ are random processes):
\begin{equation}
J_{1}(t)=\sum\limits_{i=1}^{M}\alpha_{i}E[x_{i}(t)]\label{J1}%
\end{equation}
where $\alpha_{i}$ is a weight factor for target $i$. We can now formulate a
stochastic optimization problem $\mathbf{P1}$ where the control variables are
the agent speeds and headings denoted by the vectors $\mathbf{u}%
(t)=[u_{1}(t),\dots,u_{N}(t)]$ and ${\theta}(t)=[\theta_{1}(t),\dots
,\theta_{N}(t)]$ respectively (omitting their dependence on the full system
state at $t$).
\begin{equation}
\mathbf{P1:}\qquad\qquad\quad\min\limits_{\mathbf{u(t),\theta(t)}}%
J(T)=\frac{1}{T}\int_{0}^{T}J_{1}(t)dt\quad\quad
\label{GenOptim}
\end{equation}
where $0\leq u_{j}(t)\leq1$, $0\leq\theta_{j}(t)<2\pi$, and $T$ is a given
finite mission time. This problem can be readily placed into the general
framework \eqref{GeneralJt}. In particular, the right hand side of
\eqref{GenOptim} is:%
\begin{equation}%
\begin{split}
&  \frac{1}{T}E\left[  \int_{0}^{T}\sum_{i}\int_{C(w_{i})}\frac{\alpha_{i}%
}{\pi r_{i}^{2}}x_{i}(t)dwdt\right]  \\
= &  \frac{1}{T}E\left[  \int_{0}^{T}\int_{S}\sum_{i}\frac{\alpha
_{i}\bm{1}\{w\in C(w_{i})\}}{\pi r_{i}^{2}}x_{i}(t)dwdt\right]
\end{split}
\label{J1General}%
\end{equation}
This is now in the form of the general framework in \eqref{GeneralJt} with
\begin{equation}
R(w,t)=\sum_{i}\frac{\alpha_{i}\bm{1}\{w\in C(w_{i})\}}{\pi r_{i}^{2}%
}x_{i}(t)\label{J1R}%
\end{equation}
and
\begin{equation}
P(s_{j}(t),w)=1\label{J1P}%
\end{equation}
Recalling the definition in \eqref{Rfunction_finite}, only points within the
sensing range of each target have non-zero values, while all other point value
are zero, which is the case in (\ref{J1R}) above. In addition, \eqref{J1P}
simply shows that there is no meaningful dynamic interaction between an agent
and the environment.

Problem $\mathbf{P1}$ is a finite time optimal control problem. In order to
solve this, following previous work in \cite{Khazaeni2015} we proceed with a
standard Hamiltonian analysis leading to a Two Point Boundary Value Problem
(TPBVP) \cite{bryson1975applied}. We omit this, since the details are the same
as the analysis in \cite{Khazaeni2015}. The main result of the Hamiltonian
analysis is that the optimal speed is always the maximum value, i.e.,
\begin{equation}
u_{j}^{\ast}(t)=1\label{OptimU}%
\end{equation}
Hence, we only need to calculate the optimal $\theta_{j}(t)$. This TPBVP is
computationally expensive and easily becomes intractable when problem size
grows. The ultimate solution of the TPBVP is a set of agent trajectories that
can be put in a parametric form defined by a parameter vector
$\boldsymbol{\theta}$ and then optimized over $\boldsymbol{\theta}$. If the
parametric trajectory family is broad enough, we can recover the true optimal
trajectories; otherwise, we can approximate them within some acceptable
accuracy. Moreover, adopting a parametric family of trajectories and seeking
an optimal one within it has additional benefits: it allows trajectories to be
periodic, often a desirable property, and it allows one to restrict solutions
to trajectories with desired features that the true optimal may not have,
e.g., smoothness properties to achieve physically feasible agent motion.

Parameterizing the trajectories and using gradient based optimization
methods, in light of the discussions from the previous sections, enables us to
make use of Infinitesimal Perturbation Analysis (IPA) \cite{Cassandras2010} to
carry out the trajectory optimization process. We represent each agent's
trajectory through general parametric equations
\begin{equation}%
\begin{array}
[c]{ll}%
s_{j}^{x}(t)=f_{x}(\boldsymbol{\theta}_{j},\rho_{j}(t)),\text{ \ }\quad
s_{j}^{y}(t)=f_{y}(\boldsymbol{\theta}_{j},\rho_{j}(t)) &
\end{array}
\label{param_traj}%
\end{equation}
where the function $\rho_{j}(t)$ controls the position of the agent on its
trajectory at time $t$ and $\boldsymbol{\theta}_{j}$ is a vector of parameters
controlling the shape and location of the trajectory. Let $\boldsymbol{\theta
}=[\boldsymbol{\theta}_{1},\dots,\boldsymbol{\theta}_{N}]$. We now revisit
problem $\mathbf{P1}$ in (\ref{GenOptim}):
\begin{equation}%
\begin{split}
\min\limits_{\boldsymbol{\theta}\in\Theta} J(\boldsymbol{\theta},T)= \frac
{1}{T}\int_{0}^{T}J_{1}(\boldsymbol{\theta},t)dt
\end{split}
\label{ParamOptim}%
\end{equation}
and will bring in the equations that were introduced in the previous section
in order to calculate an estimate of $\frac{dJ(\boldsymbol{\theta}%
)}{d\boldsymbol{\theta}}$ as in \eqref{Jgradient}. For this problem due to the
continuity of $x_{i}(t)$ the last two terms in \eqref{IPAsampleder} vanish.
From \eqref{J1} we have:
\begin{equation}%
\begin{split}
\frac{d}{d\boldsymbol{\theta}} \int_{\tau_{k}}^{\tau_{k+1}}\sum_{i=1}^{M}
\alpha_{i} x_{i}(\boldsymbol{\theta}, t)dt =\int_{\tau_{k}}^{\tau_{k+1}}%
\sum_{i=1}^{M} \alpha_{i} x_{i}^{\prime}(\boldsymbol{\theta},
t)dt\label{IPAsumXder}%
\end{split}
\end{equation}

In summary, the evaluation of \eqref{IPAsumXder} requires the state
derivatives $x_{i}^{\prime}(t)$ explicitly and $s_{j}^{\prime}(t)$ implicitly,
(dropping the dependence on $\boldsymbol{\theta}$ for brevity). The latter are
easily obtained for any specific choice of $f$ and $g$ in \eqref{param_traj}.
The former require a rather laborious use of
\eqref{IPAderjump},\eqref{xprime},\eqref{taukprime} which, reduces to a simple
set of state derivative dynamics as shown next.

\begin{proposition}
After an event occurrence at $t=\tau_{k}$, the state derivatives
$x_{i}^{\prime}(\tau_{k}^{+})$ with respect to the controllable parameter
$\boldsymbol{\theta}$ satisfy the following:%

\[
x_{i}^{\prime}(\tau_{k}^{+})=\left\{
\begin{array}
[c]{ll}%
0 & \text{if }e(\tau_{k})=\xi_{i}^{0}\\
x_{i}^{\prime}(\tau_{k}^{-})-\mu_{il}(t)p_{il}(\tau_{k}){\tau_{k}^{^{\prime}}}
& \text{if }e(\tau_{k})=\delta_{ij}^{+}\\
x_{i}^{\prime}(\tau_{k}^{-}) & \text{otherwise}%
\end{array}
\right.
\]
where $l\neq j$ with $p_{il}(\tau_{k})>0$ if such $l$ exists and ${\tau
_{k}^{^{\prime}}=}\frac{\partial d_{ij}(s_{j})}{\partial s_{j}}s^{\prime}%
_{j}\left(  \frac{\partial d_{ij}(s_{j})}{\partial s_{j}}\dot{s}_{j}(\tau
_{k})\right)  ^{-1}$.

\begin{proof}
The proof is omitted due to space limitations, but it is very similar to the proofs of Propositions 1-3 in \cite{Khazaeni2015rXive}.
\end{proof}
\end{proposition}

As is obvious from Proposition 1, the evaluation of $x_{i}^{\prime}(t)$ is
entirely dependent on the occurrence of events $\xi_{i}^{0}$ and $\delta
_{ij}^{+}$ in a sample realization, i.e., $\xi_{i}^{0}$ and $\delta_{ij}^{+}$
cause jumps in this derivative which carry useful information. Otherwise,
$x_{i}^{\prime}(\tau_{k}^{+})=x_{i}^{\prime}(\tau_{k}^{-})$ is in effect and
these gradients remain unchanged. However, we can easily have realizations
where no events occur in the system (specifically, events of type $\delta
_{ij}^{0}$ and $\delta_{ij}^{+}$) if the trajectory of agents in the sample
realization does not pass through any target. This lack of event excitation
results in the algorithm in \eqref{gradientmethod} to stall.

In the next section we overcome the problem of no event excitation using the
definitions in \eqref{Rfunction} and \eqref{Pfunction}. We accomplish this by
adding a new metric to the objective function that generates a non-zero
sensitivity with respect to $\boldsymbol{\theta}$.

\subsection{Event Excitation}

Our goal here is to select a function $h_{i}(\cdot )$ in \eqref{Rfunction}
with the property of \textquotedblleft spreading\textquotedblright\ the
value of $x_{i}(t)$ over all $w\in S$. We begin by determining the convex
hull produced by the targets, since the trajectories need not go outside
this convex hull. Let $\mathcal{T}=\{w_{1},w_{2},\cdots ,w_{M}\}$ be the set
of all target points. Then, the convex hull of these points is
\begin{equation}
\mathcal{C}=\bigg\{\sum_{i=1}^{M}\beta _{i}w_{i}|\sum_{i}\beta
_{i}=1,\forall i,~\beta _{i}\geq 0\bigg\}
\end{equation}%
Given that $\mathcal{C}\subset S$, we seek some $R(w,t)$ that satisfies the
following property for constants $c_{i}>0$:
\begin{equation}
\int_{\mathcal{C}}R(w,t)dw=\sum_{i=1}^{M}c_{i}x_{i}(t)  \label{Rproperty}
\end{equation}%
so that $R(w,t)$ can be viewed as a continuous density defined for all
points $w\in \mathcal{C}$ which results in a total value equivalent to a
weighted sum of the target states $x_{i}(t)$, $i=1,\ldots ,M$. In order to
select an appropriate $h(x_{i}(t),d_{i}(w))$ in \eqref{Rfunction}, we first
define $d_{i}^{+}(w)=\max (\Vert w-w_{i}\Vert ,r_{i})$ where $r_{i}$ is the
target's sensing range. We then define:
\begin{equation}
R(w,t)=\sum_{i=1}^{M}\frac{\alpha _{i}x_{i}(t)}{d_{i}^{+}(w)}  \label{J2R}
\end{equation}%
Here, we are spreading a target's reward (numerator) over all $w$ so as to
obtain the \textquotedblleft total weighted reward density" at $w$. Note
that $d_{i}^{+}(w)=\max (\Vert w-w_{i}\Vert ,r_{i})>0$ to ensure that the
target reward remains positive and fixed for points $w\in C(w_{i})$.
Moreover, following \eqref{Pfunction},
\begin{equation}
P(w,\mathbf{s}(t))=\sum_{j=1}^{N}\Vert s_{j}(t)-w\Vert ^{2}
\label{Pfunction2}
\end{equation}%
Using these definitions we introduce a new objective function metric which
is added to the objective function in \eqref{GenOptim}:
\begin{equation}
J_{2}(t)=E\Big[\int_{\mathcal{C}}P(w,\mathbf{s}(t))R(w,t)dw\Big]  \label{J2}
\end{equation}%
The expectation is a result of $P(w,\mathbf{s}(t))$ and $R(w,t)$ being
random variables defined on the same probability space as $x_{i}(t)$.

\begin{proposition}
For $R(w,t)$ in \eqref{J2R}, there exist $c_{i}>0$, $i=1,\ldots ,M$, such
that:
\begin{equation}
\int_{\mathcal{C}}R(w,t)dw=\sum_{i=1}^{M}c_{i}x_{i}(t)  \label{corequation}
\end{equation}%
\begin{proof}
  We have
  \begin{equation}
  \begin{split}
    \int_\mathcal{C}R(w,t)=\int_\mathcal{C}\sum_{i=1}^M\frac{\alpha_ix_i(t)}{d^+_i(w)}dw\\
    =\sum_{i=1}^M\alpha_i\int_\mathcal{C}\frac{x_i(t)}{d^+_i(w)}dw
    \end{split}
  \end{equation}
We now need to find the value of  $\int_\mathcal{C}\frac{x_i(t)}{d^+_i(w)}$ for each target $i$. To do this we first look at the case of one target in a $2D$ space and for now we assume $\mathcal{C}$ is just a disk with radius $\Lambda$ around the target (black circle with radius $\Lambda$ in Fig. \ref{onetargetR}). We can now calculate the above integral for this target using the polar coordinates:
\begin{equation}
\begin{split}
  &\int_\mathcal{C}\frac{x_i(t)}{d^+_i(w)}dw=\int_0^{2\pi}\int_0^{\Lambda}\frac{x_i(t)}{\max(r_i,r)}drd\theta\\
  &=\int_0^{2\pi}\int_0^{r_i}\frac{x_i(t)}{r_i}drd\theta+\int_0^{2\pi}\int_{r_i}^{\Lambda}\frac{x_i(t)}{r}drd\theta\\
  &=x_i(t)\big[2\pi\big(1+\log(\frac{\Lambda}{r_i})\big)\big]
  \end{split}
\end{equation}
\begin{wrapfigure}{r}{0.4\columnwidth}
  \centering
  \includegraphics[width=1.2in]{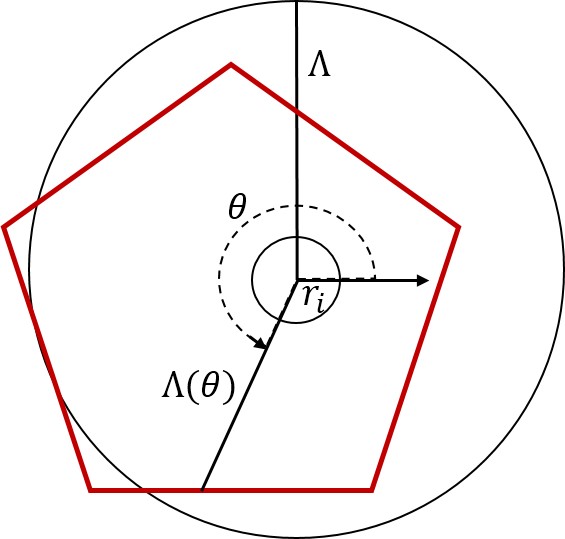}\\
  \caption{One Target $R(w,t)$ Calculation}\label{onetargetR}
\end{wrapfigure}
In our case $\mathcal{C}$ is the convex hull of all targets. We will use the same idea to calculate the $\int_\mathcal{C}\frac{x_i(t)}{d^+_i(w)}dw$ for the actual convex hull. We do this for an interior target i.e., a target inside the convex hull. Extending the same to targets on the edge is straightforward. Using the same polar coordinate for each $\theta$ we define $\Lambda(\theta)$ to be the distance of the target to the edge of $\mathcal{C}$ in the direction of $\theta$. ($\mathcal{C}$ shown by a red polygon in Fig. \ref{onetargetR}).

\begin{equation}
\begin{split}
  &\int_\mathcal{C}\frac{x_i(t)}{d^+_i(w)}dw=\int_0^{2\pi}\int_0^{\Lambda}\frac{x_i(t)}{d^+_i(r,\theta)}drd\theta\\
  &=\int_0^{2\pi}\int_0^{r_i}\frac{x_i(t)}{r_i}drd\theta+\int_0^{2\pi}\int_{r_i}^{\Lambda(\theta)}\frac{x_i(t)}{r}drd\theta\\
  &=x_i(t)\big[2\pi+\int_0^{2\pi}\log(\frac{\Lambda(\theta)}{r_i})d\theta\big]
  \end{split}
  \label{Rcalcintegral}
\end{equation}
The second part in \eqref{Rcalcintegral} has to be calculated knowing $\Lambda(\theta)$ but since we assumed the target is inside the convex hull we know $\Lambda(\theta) \ge r_i$. This means $\log(\frac{\Lambda(\theta)}{r_i})>0$ and the $x_i(t)$'s multiplier is a positive value. We can define $c_i$ in \eqref{corequation} as:
\begin{equation}
  c_i=\alpha_i\big[2\pi+\int_0^{2\pi}\log(\frac{\Lambda(\theta)}{r_i})d\theta\big]
\end{equation}
\end{proof}
\end{proposition}

The significance of $J_{2}(t)$ is that it accounts for the movement of
agents through $P(w,{\mathbf{s}}(t))$ and captures the target state values
through $R(w,t)$. Introducing this term in the objective function in the
following creates a non-zero gradient even if the agent trajectories are not
passing through any targets. We now combine the two metrics in %
\eqref{GenOptim} and \eqref{J2} and define problem $\mathbf{P2}$:
\begin{equation}
\mathbf{P2:}\qquad \min\limits_{\mathbf{u(t),\theta (t)}}J(T)=\frac{1}{T}%
\int_{0}^{T}\big[J_{1}(t)+J_{2}(t)\big]dt
\end{equation}%
In this problem, the second term is responsible for adjusting the
trajectories towards the targets by creating a potential field, while the
first term is the original performance metric which is responsible for
adjusting the trajectories so as to maximize the data collected once an
agent is within a target's sensing range. It can be easily shown that the
results in \eqref{OptimU} hold for problem $\mathbf{P2}$ as well, through
the same Hamiltonian analysis presented in \cite{Khazaeni2015}. When $%
s_{j}(t)$ follows the parametric functions in \eqref{param_traj}, the new
metric simply becomes a function of the parameter vector $\boldsymbol{\theta
}$ and we have:

\begin{equation}
\min\limits_{\boldsymbol{\theta }\in \Theta }J(\boldsymbol{\theta },T)=\frac{%
1}{T}\int_{0}^{T}\big[J_{1}(\boldsymbol{\theta },t)+J_{2}(\boldsymbol{\theta
},t)\big]dt
\label{ParamOptim2}
\end{equation}
The new objective function's derivative follows the same procedure that was
described previously. The first part's derivative can be calculated from %
\eqref{IPAsumXder}. For the second part we have:
\begin{equation}
\resizebox{0.98 \columnwidth}{!}{$ \begin{split} \frac{d}{d\bsym\theta}
&\int_{\tau_k}^{\tau_{k+1}} \int_{\mathcal{C}} P(w,\bsym
\theta,t)R(w,\bsym\theta,t)dw\\& =\int_{\tau_k}^{\tau_{k+1}}
\int_{\mathcal{C}} \Big[\frac{dP(w,\bsym \theta,t)}{d\bsym
\theta}R(w,\bsym\theta,t)+P(w,\bsym \theta,t)\frac{dR(w,\bsym\theta,t)}{d
\bsym\theta}\Big] dw \end{split}$}  \label{IPAJ2der}
\end{equation}

In the previous section, we raised the problem of no events being excited in
a sample realization, in which case the total derivative in %
\eqref{IPAsumXder} is zero and the algorithm in \eqref{gradientmethod}
stalls. Now, looking at \eqref{IPAJ2der} we can see that if no events occur
the second part in the integration which involves $\frac{dR(w,\boldsymbol{%
\theta },t)}{d\boldsymbol{\theta }}$ will be zero, since $%
\sum_{i=1}^{M}x_{i}^{\prime }(t)=0$ at all $t$. However, the first part in
the integral does not depend on the events, but calculates the sensitivity
of $P(w,\mathbf{s}(t))$ in \eqref{Pfunction2} with respect to the parameter $%
\boldsymbol{\theta }$. Note that the dependence on $\boldsymbol{\theta }$
comes through the parametric description of $\mathbf{s}(t)$ through %
\eqref{param_traj}. This term ensures that the algorithm in %
\eqref{gradientmethod} does not stall and adjusts trajectories so as to
excite the desired events.

\section{Simulation Results}

\label{numerical} 
\begin{figure}
\centering
\includegraphics[height=2.1in]{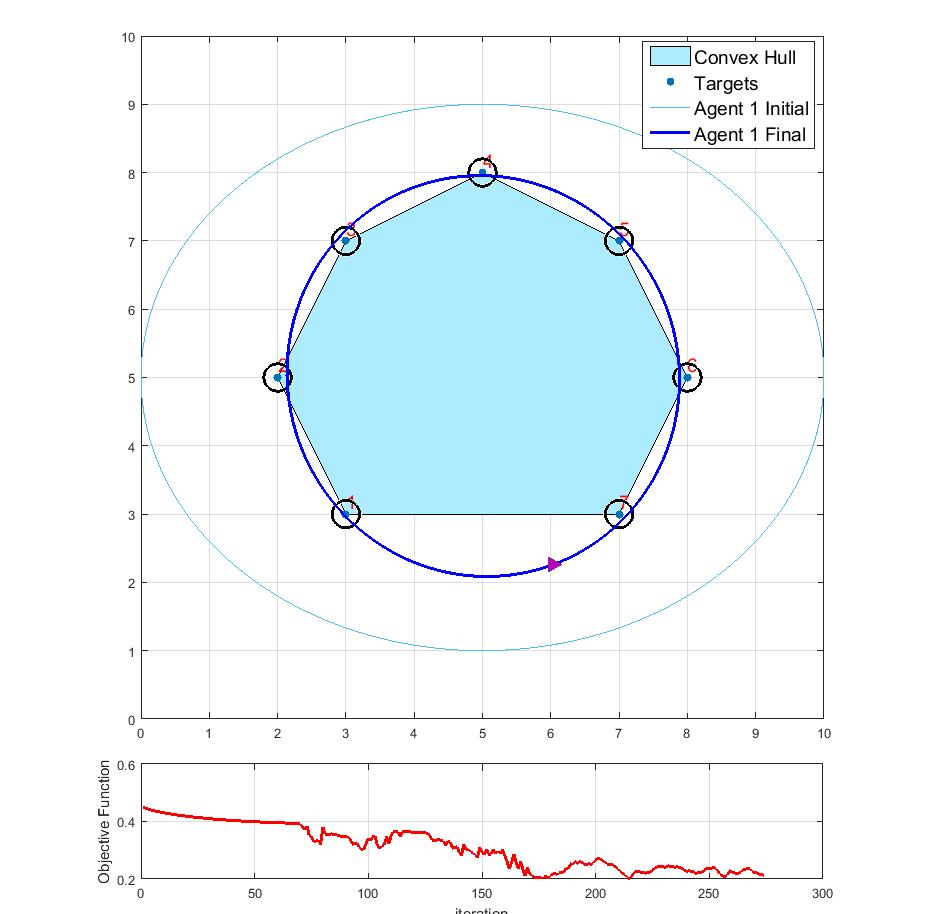}\newline
\caption{One agent and seven target scenario}
\label{oneagent}
\end{figure}

We provide some simulation results based on an elliptical
parametric description for the trajectories in \eqref{param_traj}. The
elliptical trajectory formulation is:
\begin{equation}
\begin{array}{ll}
s_{j}^{x}(t)= & A_{j}+a_{j}\cos \rho _{j}(t)\cos \phi _{j}-b_{j}\sin \rho
_{j}(t)\sin \phi _{j} \\
s_{j}^{y}(t)= & B_{j}+a_{j}\cos \rho _{j}(t)\sin \phi _{j}+b_{j}\sin \rho
_{j}(t)\cos \phi _{j}%
\end{array}%
\end{equation}%
Here, $\theta _{j}=[A_{j},B_{j},a_{j},b_{j},\phi _{j}]$ where $A_{j},B_{j}$
are the coordinates of the center, $a_{j}$ and $b_{j}$ are the major and
minor axis respectively while $\phi _{j}\in \lbrack 0,\pi )$ is the ellipse
orientation which is defined as the angle between the $x$ axis and the major
axis of the ellipse. The time-dependent parameter $\rho _{j}(t)$ is the
eccentric anomaly of the ellipse. Since an agent is moving with constant
speed of 1 on this trajectory, based on \eqref{OptimU}, we have $\dot{s}%
_{j}^{x}(t)^{2}+\dot{s}_{j}^{y}(t)^{2}=1$, which gives
\begin{equation}
\resizebox{0.99\columnwidth}{!}{$ \begin{split}
\dot\rho_j(t)=\Big[&\Big(a\sin\rho_j(t)\cos\phi_j+b_j\cos\rho_j(t)\sin\phi_j%
\Big)^2\\
&+\Big(a\sin\rho_j(t)\sin\phi_j-b_j\cos\rho_j(t)\cos\phi_j\Big)^2\Big]^{-%
\frac{1}{2}} \end{split}$}
\end{equation}%
The first case we consider is a problem with one agent and seven targets
located on a circle, as shown in Fig. \ref{oneagent}. We consider a
deterministic case with $\sigma _{i}(t)=0.5$ for all $i$. The other problem
parameters are $T=50$, $\mu _{ij}=100$, $r_{i}=0.2$ and $\alpha _{i}=1$. A
target's sensing range is denoted with solid black circles with the target
location at the center. The blue polygon indicates the convex hull produced
by the targets. The direction of motion on a trajectory is shown with the
small arrow. Starting with an initial trajectory shown in light blue, the
on-line trajectory optimization process converges to the trajectory passing
through all targets in an efficient manner (shown in dark solid blue). In
contrast, starting with this trajectory - which does not pass through any
targets - problem $\mathbf{P1}$ does not converge and the initial trajectory
remains unchanged. At the final trajectory, $J_{1}^{\ast }=0.0859$ and $%
J^{\ast }=0.2128$. Using the obvious shortest path solution, the actual
optimal value for $J_{1}$ is $0.0739$ that results from moving on the edges
of the convex hull (which allows for shorter agent travel times).

In the second case, $7$ targets are randomly distributed and two agents are
cooperatively collecting the data. The problem parameters are $\sigma
_{i}=0.5,~\mu _{ij}=10,~r_{i}=0.5,\alpha _{i}=1,~T=50$. The initial
trajectories for both agents are shown in light green and blue respectively.
We can see that both agent trajectories converge so as to cover all targets,
shown in dark green and blue ellipses. At the final trajectories, $%
J_{1}^{\ast }=0.1004$ and $J^{\ast }=0.2979$. Note that we may use these
trajectories to initialize the corresponding TPBVP, another potential
benefit of this approach. This is a much slower process which ultimately
converges to $J_{1}^{\ast }=0.0991$ and $J^{\ast }=0.2776$.
\begin{figure}
\centering
\includegraphics[height=2.1in]{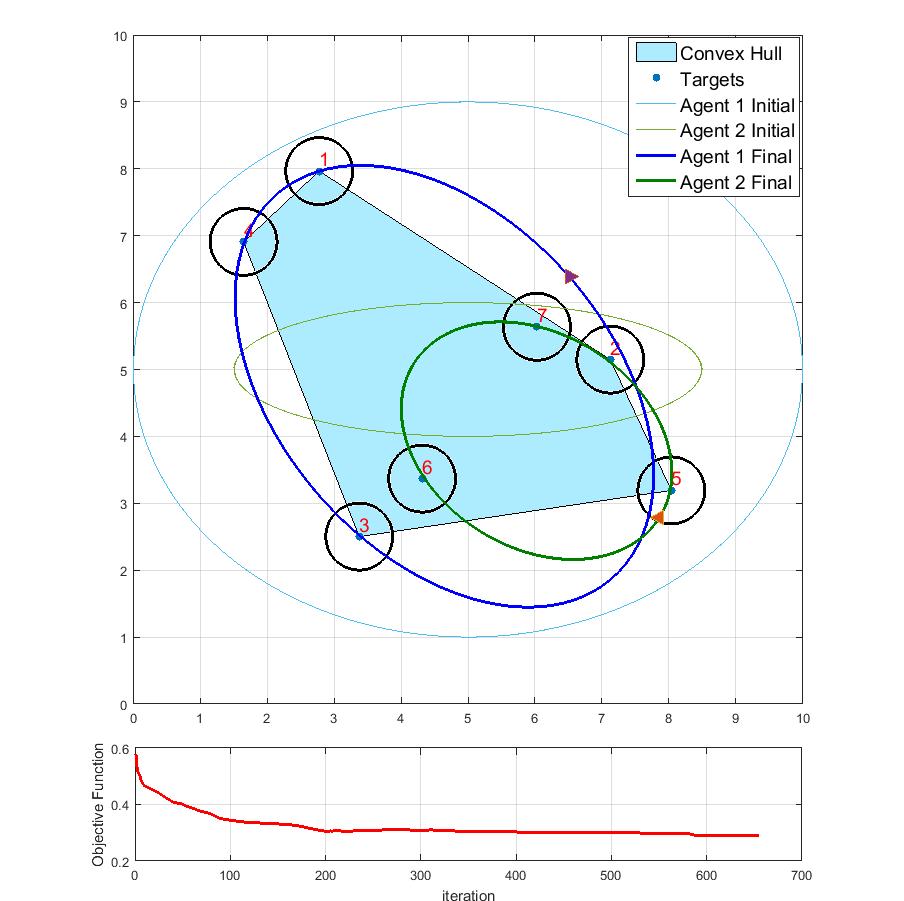}\newline
\caption{Two agent and seven targets scenario}
\label{twoagent}
\end{figure}

\section{Conclusions}

We have addressed the issue of event excitation in a class of multi-agent
systems with discrete points of interest. We proposed a new metric for such
systems that spreads the point-wise values throughout the mission space and
generates a potential field. This metric allows us to use event-driven
trajectory optimization for multi-agent systems. The methodology is applied to a class of data collection problems using the event-based IPA calculus
to estimate the objective function gradient.
\bibliographystyle{hieeetr}
\bibliography{Ref_DH}

\end{document}